

\def\ifpicture{\if01}

\input amstex
\loadbold
\documentstyle{amsppt}
\ifpicture
\input DraTex.sty
\input AlDraTex.sty

\else
\message{Please use DraTex.sty and AlDraTex.sty.
For details, see the top of this file.}
\fi
\NoBlackBoxes

\pagewidth{32pc}
\pageheight{44pc}
\magnification=\magstep1

\def\QQ{\Bbb{Q}}
\def\RR{\Bbb{R}}
\def\CC{\Bbb{C}}
\def\OO{\Cal{O}}
\def\Spec{\operatorname{Spec}}
\def\Jac{\operatorname{Jac}}
\def\rest#1#2{\left.{#1}\right\vert_{{#2}}}
\def\QED{{\unskip\nobreak\hfil\penalty50\quad\null\nobreak\hfil
{$\square$}\parfillskip0pt\finalhyphendemerits0\par\medskip}}

\topmatter
\title
Self-intersection of dualizing sheaves of
arithmetic surfaces with reducible fibers
\endtitle
\rightheadtext{}
\author Atsushi Moriwaki \endauthor
\leftheadtext{}
\address
Department of Mathematics, Faculty of Science,
Kyoto University, Kyoto, 606-01, Japan
\endaddress
\curraddr
Department of Mathematics, University of California,
Los Angeles, 405 Hilgard Avenue, Los Angeles, California 90024, USA
\endcurraddr
\email moriwaki\@math.ucla.edu \endemail
\date February, 1995 \enddate
\abstract
Let $K$ be an algebraic number field and $O_K$ the ring of integers of $K$.
Let $f : X \to \Spec(O_K)$ be a stable arithmetic surface
over $O_K$ of genus $g \geq 2$.
In this short note, we will prove that if $f$ has a reducible geometric fiber,
then
$(\omega_{X/O_K}^{Ar} \cdot \omega_{X/O_K}^{Ar} ) \geq
\log 2/6(g-1)$.
\endabstract
\endtopmatter

\document

Let $K$ be a number field and $O_K$ the ring of integers of $K$.
Let $f : X \to \Spec(O_K)$ be a stable arithmetic surface of genus $g \geq 2$.
One of interesting problems on arithmetic surfaces is a question whether
$(\omega_{X/O_K}^{Ar} \cdot \omega_{X/O_K}^{Ar} ) > 0$.
This question is closely related to Bogomolov conjecture, which claims that,
for a curve $C$ over a number field and an embedding $j : C \to \Jac(C)$ into
the Jacobian $\Jac(C)$, the image $j(C)$ in $\Jac(C)$ is discrete in terms of
the N\'{e}ron-Tate height of $\Jac(C)$ (cf. \cite{Sz} and \cite{Zh2}).
Currently, the positivity of $(\omega_{X/O_K}^{Ar} \cdot \omega_{X/O_K}^{Ar} )$
is known in the following cases:
\roster
\item "(1)" $f$ has a bad reduction. (\cite{Zh1} and \cite{Zh2})

\item "(2)" $f$ has a reducible geometric fiber, or $\Jac(X_K)$
has a complex multiplication. (\cite{Bu})

\item "(3)" $\operatorname{End}(\Jac(X_K))_{\RR}$ is not
isomorphic to $\RR$, $\CC$, or the quaternion division algebra $\Bbb{D}$.
(\cite{Zh3})
\endroster
In this short note, we would like to give an effective lower bound
of $(\omega_{X/O_K}^{Ar} \cdot \omega_{X/O_K}^{Ar} )$ under the assumption that
$f$ has a reducible geometric fiber, namely,

\proclaim{Theorem 1}
Let $K$ be an algebraic number field and $O_K$ the ring of integers of $K$.
Let $f : X \to \Spec(O_K)$ be a stable arithmetic surface
over $O_K$ of genus $g \geq 2$.
Assume that, for $P_1, \ldots, P_n \in \Spec(O_K)$,
geometric fibers $X_{\overline{P}_1}, \ldots, X_{\overline{P}_n}$ of $X$
at $P_1, \ldots, P_n$ are reducible.
Then, we have
$$
\left(\omega_{X/O_K}^{Ar} \cdot \omega_{X/O_K}^{Ar} \right) \geq
\sum_{i=1}^n \frac{\log \#(O_K/P_i) }{6(g-1)}.
$$
In particular, $\left(\omega_{X/O_K}^{Ar} \cdot \omega_{X/O_K}^{Ar} \right)
\geq
\log 2/6(g-1)$
if $f$ has a reducible geometric fiber.
\endproclaim

\noindent
As compared with methods used in \cite{Bu}, \cite{Zh1}, \cite{Zh2} and
\cite{Zh3}, our method is very elementary.
Besides it, we can give the exact positive lower bound
$\log 2 / 6(g-1)$.

\bigskip
Let $h : Y \to \Spec(O_K)$ be a regular semistable arithmetic surface over
$O_K$
of genus $g \geq 2$. Then, $\omega_{Y/O_K}$ is $h$-nef, i.e.,
$\deg(\rest{\omega_{Y/O_K}}{C}) \geq 0$ for all vertical curves $C$ on $Y$.
Let $P \in \Spec(O_K)$ and $Y_P$ the fiber of $Y$ at $P$.
To distinguish Arakelov intersection and
geometric intersection in the fiber $Y_P$, we will use
a symbol $(\ ,\ )$ for Arakelov intersection and
$\langle\ ,\ \rangle$ for geometric intersection.
$(\ ,\ )$ and $\langle\ ,\ \rangle$ are related by a formula
$$
(\ ,\ ) = \langle\ ,\ \rangle\log\#(O_K/P).
$$
Let $Y_{\overline{P}} = Y_P \otimes \overline{(O_K/P)}$
be the geometric fiber at $P$ and
$C$ an irreducible component of $Y_{\overline{P}}$.
We say $C$ is an {\it essential irreducible component}
if $\langle \omega_{Y/O_K} \cdot C \rangle > 0$.
Here we assume that $Y_{\overline{P}}$ is reducible.
We set
$$
\beta(P) = \min \left.\left\{
- \frac{\langle \omega_{Y/O_K} \cdot C \rangle^2}
{\langle C \cdot C \rangle} \right|
\text{$C$ is an essential irreducible component of $Y_{\overline{P}}$.}
\right\}.
$$

\proclaim{Lemma 2} $\beta(P) \geq 1/3$.
\endproclaim

\demo{Proof}
Let $C$ be an irreducible component of $Y_{\overline{P}}$ with
$\langle \omega_{Y/O_K} \cdot C \rangle > 0$.
We need to prove
$$
\frac{\langle \omega_{Y/O_K} \cdot C \rangle^2}
{\langle C \cdot C \rangle} \leq -\frac{1}{3}.
$$
By adjunction formula, if $q = h^1(\OO_C)$ is an arithmetic genus of $C$, then
$$
\langle \omega_{Y/O_K} \cdot C \rangle + \langle C \cdot C \rangle
= 2q - 2.
$$
Thus,
$$
\frac{\langle \omega_{Y/O_K} \cdot C \rangle^2}
{\langle C \cdot C \rangle} = \langle C \cdot C \rangle - 4(q-1)
+ \frac{4(q-1)^2}{\langle C \cdot C \rangle}.
$$
Therefore, since $\langle C \cdot C \rangle \leq -1$,
if $q \geq 1$, our assertion is trivial.

Here we assume $q = 0$. Then, the previous formula is
$$
\frac{\langle \omega_{Y/O_K} \cdot C \rangle^2}
{\langle C \cdot C \rangle} = \langle C \cdot C \rangle + 4 +
\frac{4}{\langle C \cdot C \rangle}.
$$
Since $\langle \omega_{Y/O_K} \cdot C \rangle > 0$,
we have $\langle C \cdot C \rangle \leq -3$.
Therefore, we get our lemma.
\QED
\enddemo

Let us consider the following conditions for $Y_P$:
\roster
\item "(C.1)"
Every irreducible component of $Y_P$ is geometrically irreducible and
every node of the geometric fiber $Y_{\overline{P}}$ is defined over
$O_K/P$.

\item "(C.2)"
For any point closed $x$ of the fiber $Y_{{P}}$,
there is an irreducible component $C$ of $Y_P$ such that
$\langle \omega_{Y/O_K} \cdot C \rangle > 0$ and $x \not\in C$.
\endroster
Note that (C.2) holds if (C.1) holds
and number of essential irreducible components
is greater than or equal to $3$.

\proclaim{Lemma 3}
Let $\Gamma$ be a horizontal prime divisor on $Y$, and
$\rho : \widetilde{\Gamma} = \Spec(O_{K'}) \to \Gamma$
a normalization of $\Gamma$.
Let $Y' \to Y \times_{\Spec(O_K)} \widetilde{\Gamma}$
be a minimal resolution of $Y \times_{\Spec(O_K)} \widetilde{\Gamma}$.
We set induced morphisms as follows:
$$
\CD
Y @<{\rho'}<< Y' \\
@V{h}VV @VV{h'}V \\
\Spec(O_K) @<<{\rho}< \widetilde{\Gamma} = \Spec(O_{K'})
\endCD
$$
Let $\Gamma'$ be a section of $g'$ with $\rho'(\Gamma') = \Gamma$.
Under the assumption \rom{(C.1)} and \rom{(C.2)} for $Y_P$,
there is a vertical $\QQ$-divisor $E_P$ on $Y'$ with following properties:
\roster
\item "3.1." $\rho'(E_P) \subseteq Y_P$.

\item "3.2." $\left( (\omega^{Ar}_{Y'/O_{K'}} - (2g-2) \overline{\Gamma}' -
E_P)
\cdot E_P \right)= 0$.

\item "3.3." ${\displaystyle (E_P \cdot E_P) \leq
- \beta(P) \log\#(O_K/P) [K' : K]}$.
\endroster
\endproclaim

\demo{Proof}
First of all, it is well known that
$\omega_{Y'/O_{K'}} = {\rho'}^*(\omega_{Y/O_K})$.
We set
$$
PO_{K'} = Q_1^{e_1} \cdots Q_t^{e_t}.
$$
For each $i$ ($1 \leq i \leq t$),
we would like to find $\QQ$-divisors $E_i$ on $Y$ with
the following properties:
\roster
\item "(1)" $E_i$ is a scalar of an essential irreducible component
$C$ of $Y_P$ by a rational number, i.e.,
$E_i = a C$ for some $a \in \QQ$.

\item "(2)" $\langle E_i \cdot E_i \rangle \leq -\beta(P)$.

\item "(3)" If we set $E'_i = \rest{{\rho'}^*(E_i)}{Y'_{Q_i}}$, then
$\left( (\omega^{Ar}_{Y'/O_{K'}} - (2g-2) \overline{\Gamma}' - E'_i)
\cdot E'_i \right)= 0.$
\endroster
Assuming the existence of $E_i$'s as above, we proceed to our argument.
We set
$$
E_P = E'_1 + \cdots + E'_t.
$$
Then, 3.1 and 3.2 are obvious by (1) and (3).
Let us try to evaluate $\langle E'_i \cdot E'_i \rangle$.
Let $C'$ be an irreducible component of $Y'_{Q_i}$ such that
$\rho'(C') = C$.
Then,
$$
E'_i = e_i a C' + (\text{exceptional divisors}).
$$
Thus, if we set $\rho'_*(C') = f C$, then
$$
\align
\langle E'_i \cdot E'_i \rangle
& = e_i  a \langle \rest{{\rho'}^*(E_i)}{Y'_{Q_i}} \cdot C' \rangle
= e_i a \langle E_i \cdot \rho'_*(C') \rangle \\
& = e_i a f \langle E_i \cdot C \rangle
= e_i f \langle E_i \cdot E_i \rangle
\leq -\beta(P) e_i .
\endalign
$$
Therefore,
$$
\align
(E_P \cdot E_P) & =
\sum_{i=1}^t (E'_i \cdot E'_i)  =
\sum_{i=1}^t  \langle E'_i \cdot E'_i \rangle \log\#(O_{K'}/Q_i) \\
& \leq \sum_{i=1}^t -\beta(P) e_i \log\#(O_{K'}/Q_i) \\
& = -\beta(P) \log\#(O_K/P) \sum_{i=1}^t e_i [O_{K'}/Q_i :  O_K/P] \\
& = -\beta(P) \log\#(O_K/P)  [K' : K].
\endalign
$$
Thus $E_P$ satisfies 3.3.

We need to construct $E_i$ satisfying (1), (2) and (3).
By our assumptions (C.1) and (C.2),
there is an essential irreducible component $C$ of $Y_P$
such that $\rho'(Y'_{Q_i} \cap \Gamma') \not\in C$.
We set
$$
E_i = \frac{\langle \omega_{Y/O_K} \cdot C \rangle}
{\langle C \cdot C \rangle} C.
$$
Then, by the definition of $\beta(P)$,
$\langle E_i \cdot E_i \rangle \leq -\beta(P)$.
Let $C'$ be an irreducible component of $Y'_{Q_i}$ such that
$\rho'(C') = C$.
Then, we can write
$$
E'_i = \rest{{\rho'}^*(E_i)}{Y'_{Q_i}} = b C' + F,
$$
where $b$ is a rational number and
$F$ is an exceptional divisor of $\rho'$.
By our choice of $C$, $C' \cap \Gamma' = \emptyset$ and
$F \cap \Gamma' = \emptyset$.
Therefore, if we set $\rho'_*(C') = fC$, then
$$
\align
\left\langle (\omega_{Y'/O_{K'}} - (2g-2) \overline{\Gamma}' - E'_i)
\cdot E'_i \right\rangle & =
\left\langle (\omega_{Y'/O_{K'}} - (2g-2) \overline{\Gamma}' - E'_i)
\cdot b C' \right\rangle \\
& =
\left\langle (\omega_{Y'/O_{K'}} - E'_i)
\cdot b C' \right\rangle \\
& = b \langle \omega_{Y/O_K} \cdot \rho'_*(C') \rangle -
    b \langle E_i \cdot \rho'_*(C') \rangle \\
& = b f \left(
\langle \omega_{Y/O_K} \cdot C \rangle -
    \langle E_i \cdot C \rangle \right)  = 0.
\endalign
$$
which completes the proof of Lemma~3.
\QED
\enddemo

Using Lemma~3, we get

\proclaim{Lemma 4}
Let $P_1, \ldots, P_n \in \Spec(O_K)$ such that geometric
fibers $Y_{\overline{P}_1}, \ldots, Y_{\overline{P}_n}$ at $P_1, \ldots, P_n$
are reducible. If $Y_{P_1}, \ldots, Y_{P_n}$ satisfy the conditions \rom{(C.1)}
and
\rom{(C.2)}, then,
for any prime divisor $\Gamma$
\rom{(}including a infinite fiber\rom{)} on $Y$, we have
$$
\left[ \left(\omega_{Y/O_K}^{Ar} \cdot \omega_{Y/O_K}^{Ar} \right) +
\sum_{i=1}^n \beta(P_i) \log\#(O_K/P_i) \right] \deg(\Gamma_K) \leq
2g(2g-2) \left(\omega_{Y/O_K}^{Ar} \cdot \overline{\Gamma}\right).
$$
\endproclaim

\demo{Proof}
Since $\omega_{Y/O_K}$ is $h$-nef,
if $\Gamma$ is an irreducible component of a fiber, or
an infinite fiber, our assertion is trivial.
So we assume that $\Gamma$ is horizontal.
We use the same notation as in Lemma~3.
By virtue of Lemma~3, for each $P_i$, there is a vertical $\QQ$-divisor
$E_{P_i}$
on $Y'$ with the following properties:
\roster
\item "(a)" $\rho'(E_{P_i}) \subseteq Y_{P_i}$.

\item "(b)" $\left( (\omega^{Ar}_{Y'/O_{K'}} - (2g-2) \overline{\Gamma}' -
E_{P_i})
\cdot E_{P_i} \right)= 0$.

\item "(c)" ${\displaystyle (E_{P_i} \cdot E_{P_i}) \leq
- \beta(P_i) \log\#(O_K/P_i) [K' : K]}$.
\endroster
We set $E = E_{P_1} + \cdots + E_{P_n}$.
Then, by (a), (b), (c) and Hodge index theorem (\cite{Fa, Theorem~4}),
$$
\align
\left( \omega^{Ar}_{Y'/O_{K'}} - (2g-2) \overline{\Gamma}' \right)^2 &=
\left( \omega^{Ar}_{Y'/O_{K'}} - (2g-2) \overline{\Gamma}' - E + E \right)^2 \\
& =
\left( \omega^{Ar}_{Y'/O_{K'}} - (2g-2) \overline{\Gamma}' - E \right)^2
+ (E \cdot E) \\
& \leq \sum_{i=1}^n (E_{P_i} \cdot E_{P_i}) \\
& \leq \sum_{i=1}^n -\beta(P_i) \log\#(O_K/P_i) [K' : K].
\endalign
$$
It follows from adjunction formula that
$$
(\omega^{Ar}_{Y'/O_{K'}} \cdot \omega^{Ar}_{Y'/O_{K'}})
- 2g(2g-2) (\omega^{Ar}_{Y'/O_{K'}} \cdot \overline{\Gamma}')
\leq -\sum_{i=1}^n \beta(P_i) \log\#(O_K/P_i)  [K' : K].
$$
Hence, we obtain
$$
\left[ \left(\omega_{Y/O_K}^{Ar} \cdot \omega_{Y/O_K}^{Ar} \right) +
\sum_{i=1}^n \beta(P_i) \log\#(O_K/P_i) \right] \deg(\Gamma_K) \leq
2g(2g-2) \left(\omega_{Y/O_K}^{Ar} \cdot \overline{\Gamma}\right).
$$
\QED
\enddemo

We need the following technical fact of algebra.

\proclaim{Lemma 5}
Let $K$ be an algebraic number field and $O_K$ the ring of integers of $K$.
Let $P_1, \ldots, P_r$ be non-zero prime ideals of $O_K$ and
$n$ a positive integer. Then, there is a finite extension field $K'$ of $K$
with the following properties:
\roster
\item "(1)" $[K' : K] = n$.

\item "(2)" $Q_i = P_i O_{K'}$ is a prime ideal for every $i$.

\item "(3)" $[O_{K'}/Q_i : O_K/P_i] = n$ for every $i$.
\endroster
\endproclaim

\demo{Proof}
Let $(O_K)_{P_i}$ be the localization of $O_K$ at $P_i$ and
$F_i \in (O_K)_{P_i}[X]$ such that $F_i$ is monic and of degree $n$, and
the class $\bar{F}_i$ in $(O_K)_{P_i}/P_i(O_K)_{P_i}[X]$
($= O_K/P_i[X]$) is irreducible.
We set
$$
F_i = X^n + a_{i1}X^{n-1} + \cdots + a_{in}.
$$
By approximation theorem (cf. \cite{La, XII, Theorem~1.2}),
for every $1 \leq j \leq n$,  there is $a_j \in K$ such that
$$
a_j - a_{ij}  \in P_i(O_K)_{P_i}
$$
for all $1 \leq i \leq r$.
We set
$$
F = X^n + a_1 X^{n-1} + \cdots + a_n.
$$
Then, $F \in (O_K)_{P_i}[X]$ and $F \equiv F_i \pmod{P_i(O_K)_{P_i}[X]}$
for all $1 \leq i \leq r$.
Since $\bar{F} = \bar{F}_i$ in $(O_K)_{P_i}/P_i(O_K)_{P_i}[X]$
is irreducible, $F$ is irreducible over $K$ by
Gauss Lemma (cf. \cite{La, IV, Corollary~2.2}).
Let $\alpha$ be a root of $F$ and $K' = K(\alpha)$. Then, $[K' : K] = n$.
Let $Q_i$ be a prime ideal of $O_{K'}$ lying
over $P_i$. Since $\alpha$ is integral over $(O_K)_{P_i}$,
$\alpha \in (O_{K'})_{Q_i}$. Let
$\bar{\alpha}$ be the class of $\alpha$ in
$O_{K'}/Q_i = (O_{K'})_{Q_i}/Q_i(O_{K'})_{Q_i}$.
Then, since $\bar{\alpha}$ is a root of $\bar{F}_i$ and $\bar{F}_i$ is
irreducible,
$[O_{K'}/Q_i : O_K/P_i] \geq n$.
On the other hand, by \cite{La, XII, Corollary~6.3},
$[O_{K'}/Q_i : O_K/P_i] \leq n$ and equality holds if and only if
$Q_i = P_i O_{K'}$.
Thus, we get our lemma.
\QED
\enddemo

\bigskip
Let us start the proof of Theorem~1.
Actually, we will prove the following.
Indeed, Theorem~6 implies Theorem~1 if we consider
a minimal resolution $Y \to X$.

\proclaim{Theorem 6}
Let $h : Y \to \Spec(O_K)$ be a regular semistable arithmetic surface
of genus $g \geq 2$. Let $P_1, \ldots, P_n \in \Spec(O_K)$
such that, for every $1 \leq i \leq n$,
a geometric fiber $Y_{\overline{P}_i}$ at $P_i$ has  at least two
essential irreducible components. Then, we have the following:
\roster
\item "(A)" For all effective Arakelov divisors $D$ on $Y$,
$$
\left[ \left(\omega_{Y/O_K}^{Ar} \cdot \omega_{Y/O_K}^{Ar} \right) +
\sum_{i=1}^n \frac{\log\#(O_K/P_i)}{6} \right] \deg(D_K) \leq
2g(2g-2) \left(\omega_{Y/O_K}^{Ar} \cdot D \right).
$$

\item "(B)" ${\displaystyle
\left(\omega_{Y/O_K}^{Ar} \cdot \omega_{Y/O_K}^{Ar} \right) \geq
\sum_{i=1}^n \frac{\log \#(O_K/P_i) }{6(g-1)}}$.
\endroster
\endproclaim

\demo{Proof}
First, let us consider (A).
Here we claim:

\proclaim{Claim 7}
We may assume that $Y_{P_i}$ satisfies \rom{(C.1)} for every $1 \leq i \leq n$.
\endproclaim

Let $k_i$ be a finite extension field of $O_K/P_i$
such that every irreducible component of $Y_{P_i} \otimes k_i$
is geometrically irreducible and
every node of $Y_{P_i} \otimes k_i$ is a $k_i$-rational point.
Let $l$ be a positive integer such that
$[k_i : O_K/P_i] \leq l$ for all $1 \leq i \leq n$.
Let $P_{n+1}, \ldots, P_r \in \Spec(O_K)$ be other critical values of $h$.
Then, by Lemma~5, there is a finite extension field $E$ of $K$ such that
$[E : K] = l$, $R_i = P_i O_{E}$ is a prime ideal and
$[O_{E}/R_i : O_K/P_i] = l$ for every $1 \leq i \leq r$.
Let us consider  $Z = Y \times_{\Spec(O_K)} \Spec(O_{E})$.
Then, $Z$ is regular and semistable, and for $1 \leq i \leq n$,
a fiber $Z_{R_i}$ at $R_i$ satisfies (C.1).
Thus, we get
$$
\left[ \left(\omega_{Z/O_{E}}^{Ar} \cdot \omega_{Z/O_{E}}^{Ar} \right) +
\sum_{i=1}^n \frac{\log\#(O_E/R_i)}{6} \right] \deg(\rho^*(D)_{E})
\leq 2g(2g-2) \left(\omega_{Z/O_{E}}^{Ar} \cdot \rho^*(D) \right)
$$
for all effective Arakelov divisors $D$ on $Y$,
where $\rho : Z \to Y$ is the induced morphism.
Here, $\omega_{Z/O_{E}}^{Ar} = \rho^*(\omega_{Y/O_{K}}^{Ar})$,
and $\log\#(O_E/R_i) = \log\#(O_K/P_i)[E : K]$.
Thus, we obtain our claim.

\medskip
Moreover, in the same way as in Claim~7, we may assume that
if we take a blowing-up at a node of $Y_{P_i}$, then
new nodes are also $(O_K/P_i)$-rational points.

Let $t_{P_i}$ be a local parameter of $P_i(O_K)_{P_i}$.
We set $K' = K(\sqrt{\prod_{i=1}^n t_{P_i}})$.
Then, there is a unique prime ideal $Q_i$ of $O_{K'}$ such that
$Q_i \cap O_K = P_i$ and the ramification index between $P_i$ and $Q_i$
is equal to $2$.
Let $Y_1 = Y \times_{\Spec(O_K)} \Spec(O_{K'})$ and
$Y' \to Y_1$ a minimal resolution of $Y_1$.
If number of essential irreducible components of $Y_{P_i}$ is
greater than or equal to $3$,
then $Y_{P_i}$ satisfies the condition (C.2). So does $Y'_{Q_i}$.
We assume that
number of essential irreducible components of $Y_{P_i}$ is equal to $2$.
Let $C_1$ and $C_2$ be essential irreducible components of $Y_{P_i}$.
Let $C'_1$ and $C'_2$ be irreducible components of a fiber $(Y_1)_{Q_i}$
corresponding to
$C_1$ and $C_2$. Then, $Y_1$ has $A_1$ type singularities at
$C'_1 \cap C'_2$. Therefore,
corresponding irreducible components $C''_1$ and $C''_2$ in $Y'_{Q_i}$
has no common points. (Note that every node of $Y'_{Q_i}$ is defined
over $O_{K'}/Q_i = O_K/P_i$.)
Thus, $Y'_{Q_i}$ satisfies the condition (C.2).
\ifpicture
\message{Making picture. It takes a long time...}
\par\medskip
\Draw
\MoveTo(0,0)
\DrawRect(50,75)
\MoveTo(90,0)
\DrawRect(50,75)
\MoveTo(180,0)
\DrawRect(70,75)
\MoveTo(60,37.5) \Node(A)(----)
\MoveTo(80,37.5) \Node(B)(----)
\MoveTo(150,37.5) \Node(C)(----)
\MoveTo(170,37.5) \Node(D)(----)
\ArrowHeads(1) \Edge(B,A) \Edge(D,C)
\MoveTo(25,-10) \Node(E)(--$Y$--)
\MoveTo(115,-10) \Node(F)(--$Y_1$--)
\MoveTo(215,-10) \Node(G)(--$Y'$--)
\MoveTo(5,37.5) \DrawOvalArc(40,35)(-90,90)
\MoveTo(45,37.5) \DrawOvalArc(40,35)(90,270)
\MoveTo(95,37.5) \DrawOvalArc(40,35)(-90,90)
\MoveTo(135,37.5) \DrawOvalArc(40,35)(90,270)
\MoveTo(182,2) \MarkLoc(1) \Move(10,10) \MarkLoc(2)
\MoveTo(182,73) \MarkLoc(3) \Move(10,-10) \MarkLoc(4)
\Curve(1,2,4,3)
\MoveTo(248,2) \MarkLoc(5) \Move(-10,10) \MarkLoc(6)
\MoveTo(248,73) \MarkLoc(7) \Move(-10,-10) \MarkLoc(8)
\Curve(5,6,8,7)
\MoveTo(12.5,37.5) \Node(H)(--$C_1$--)
\MoveTo(37.5,37.5) \Node(I)(--$C_2$--)
\MoveTo(102.5,37.5) \Node(J)(--$C'_1$--)
\MoveTo(127.5,37.5) \Node(K)(--$C'_2$--)
\MoveTo(197,37.5) \Node(L)(--$C''_2$--)
\MoveTo(233,37.5) \Node(M)(--$C''_1$--)
\MoveTo(115,7.5) \Node(N)(--$\bullet$--)
\MoveTo(115,67.5) \Node(O)(--$\bullet$--)
\MoveTo(115,15) \Node(P)(--$A_1$--)
\MoveTo(115,60) \Node(Q)(--$A_1$--)
\MoveTo(182,7.5) \Line(66,0)
\MoveTo(182,67.5) \Line(66,0)
\MoveTo(215,15) \Node(R)(--(-2)-curve--)
\MoveTo(215,60) \Node(S)(--(-2)-curve--)
\EndDraw
\medskip\par\noindent
\else
\fi
Hence we can apply Lemma~4 for $Y'$.
Namely if $\mu : Y' \to Y$ is the induced morphism, then,
$$
\multline
\left[ \left(\omega_{Y'/O_{K'}}^{Ar} \cdot \omega_{Y'/O_{K'}}^{Ar} \right) +
\sum_{i=1}^n \beta(Q_i) \log\#(O_K/Q_i) \right] \deg(\mu^*(D)_{K'}) \\
\leq 2g(2g-2) \left(\omega_{Y'/O_{K'}}^{Ar} \cdot \mu^*(D) \right).
\endmultline
$$
Thus, by Lemma~2, we get (A) because
$\omega_{Y'/O_{K'}}^{Ar} = \mu^*(\omega_{Y/O_{K}}^{Ar})$ and
$[K' : K] = 2$.

\medskip
Next let us consider (B).
Let $F$ be an infinite fiber and $r$ a real number with
$$
(\omega^{Ar}_{Y/O_K} - r F \cdot \omega^{Ar}_{Y/O_K} - r F) > 0,
$$
i.e. $r < \frac{(\omega^{Ar}_{Y/O_K} \cdot \omega^{Ar}_{Y/O_K})}
{2 (\omega^{Ar}_{Y/O_K} \cdot F) }$.
Then, by \cite{Fa, Corollary of Theorem~3},
there is an effective Arakelov divisor $E$ on $Y$
such that $E$ is arithmetically linearly equivalent to
$n(\omega^{Ar}_{Y/O_K} - r F)$ for sufficiently large $n$.
Thus, by (A), we have
$$
\left(\omega_{Y/O_K}^{Ar} \cdot \omega_{Y/O_K}^{Ar} \right) +
\sum_{i=1}^n \frac{\log\#(O_K/P_i)}{6} \leq
2g \left(\omega_{Y/O_K}^{Ar} \cdot
\omega^{Ar}_{Y/O_K} - r F \right).
$$
Taking a limit
$r \to \frac{(\omega^{Ar}_{Y/O_K} \cdot \omega^{Ar}_{Y/O_K})}
{2 (\omega^{Ar}_{Y/O_K} \cdot F) }$,
we get (B).
\QED
\enddemo


\enddocument